\documentclass{ifacconf}

\usepackage{graphicx}      
\usepackage[round,authoryear,sort]{natbib}        
\usepackage{graphics} 
\usepackage{epsfig} 
\usepackage{amsmath} 
\usepackage{amssymb}  
\usepackage{algorithm}
\usepackage{algpseudocode}
\usepackage{mathtools}

\usepackage{caption}
\usepackage{subcaption}
\newcommand{\set}[1]{\left\{#1\right\}}
\newcommand{\bs}{\boldsymbol}

\newcommand{\norm}[1]{\left\Vert#1\right\Vert}
\begin{document}
\begin{frontmatter}

\title{Numerical discrete-time implementation of continuous-time linear-quadratic model predictive control} 


\author{Zhanhao Zhang}, 
\author{Anders H.D. Christensen},
\author{Steen Hørsholt}, 
\author{John Bagterp Jørgensen}

\address{Department of Applied Mathematics and Computer Science, Technical University of Denmark, DK-2800 Kgs. Lyngby, Denmark (e-mail: jbjo@dtu.dk).}

\begin{abstract}                

This study introduces the design, discretization, and discrete-time implementation of continuous-time linear-quadratic model predictive control (CT-LMPC). The control model for CT-LMPC is parameterized using transfer functions with time delays, which are decomposed into deterministic and stochastic parts to facilitate appropriate control and filtering algorithms. The CT-LMPC addresses continuous-time (CT), time-delayed, linear-quadratic optimal control problems (LQ-OCPs). By considering piecewise constant inputs, we propose a numerical discretization method for LQ-OCPs and demonstrate how their discrete-time (DT) equivalents can be transformed into a standard quadratic program. The performance of CT-LMPC is evaluated against that of conventional DT-LMPC. Numerical experiments reveal that, with fixed tuning parameters, CT-LMPC exhibits superior closed-loop performance compared to DT-LMPC as the sampling time increases.

\end{abstract}


\end{frontmatter}
\section{Introduction}
\label{sec:Introduction}

Consider the following continuous-time (CT) linear-quadratic optimal control problem (LQ-OCP) with piecewise constant inputs 
\begin{subequations}
\label{eq:continuousLQOCP}
\begin{alignat}{3}
    & \min \: \: && \phi=\int_{t_0}^{t_0+T} l(\tilde z(t)) dt \\
& s.t. && x(t_0) = \hat x_0, \\
& && u(t) = u_k,  &&t_k \leq t < t_{k+1}, \, k \in \mathcal{N}, \label{eq:ZOHInput}\\
& && \dot x(t) = A_c x(t) + B_c u(t), \,  && t_0 \leq t < t_0+T, \\
& && z(t) = C_c x(t) + D_c u(t), && t_0 \leq t < t_0+T, \\
& && \bar z(t) = \bar z_k,  &&t_k \leq t < t_{k+1}, \, k \in \mathcal{N}, \label{eq:ZOHOutputRef} \\
& && \tilde z(t) = z(t) - \bar z(t), && t_0 \leq t < t_0+T,
\end{alignat}
\end{subequations}
with the stage cost function 
\begin{equation}
\begin{split}
    l_c(\tilde z(t)) &= \frac{1}{2} \norm{ W_z \tilde z(t)  }_2^2
    = \frac{1}{2} \left( \tilde z(t) ' Q_{c} \tilde z(t) \right), 
\end{split}
\label{eq:det-stagecost}
\end{equation}
where $\bar z$ is the reference and is assumed to be piecewise constant. $Q_{c} = W_z' W_z$ is semi-positive definite weight matrix. $T=NT_s$ is the control horizon and $\mathcal{N}=0,1,\ldots,N-1$.

The corresponding discrete-time (DT) equivalent of~\eqref{eq:continuousLQOCP} is
\begin{subequations}
\label{eq:discreteLQOCP}
\begin{alignat}{5}
& \min_{ x,u} \quad &&\phi = \sum_{k\in \mathcal{N}} l_k(x_k,u_k)  \\
& s.t. && x_0 = \hat x_0, \\
& && x_{k+1} = A x_k + B u_k, \quad && k \in \mathcal{N},
\end{alignat}    
\end{subequations}
with the stage costs
\begin{equation}
    l_k(x_k,u_k) = \frac{1}{2} \begin{bmatrix} x_k \\ u_k \end{bmatrix}' Q \begin{bmatrix} x_k \\ u_k \end{bmatrix} + q_k' \begin{bmatrix} x_k \\ u_k \end{bmatrix} + \rho_k, \quad k \in \mathcal{N},
\label{eq:deterministic-stageCost}
\end{equation}
where $Q$ is a symmetric positive semi-definite matrix and 
\begin{equation}
\label{eq:analyticQandq}
    Q = \begin{bmatrix} Q_{xx} & Q_{xu} \\ Q_{ux} & Q_{uu} \end{bmatrix},  \quad 
    q_k = M \bar z_k, \quad \rho_k = l_c(\bar z_k) T_s. 
\end{equation}
Linear-quadratic optimal control problems (LQ-OCPs) are considered one of the fundamental problems in optimal control theory due to their simplicity and mathematical tractability. There is abundant research on the solution and discretization methods of optimal control problems described in~\eqref{eq:continuousLQOCP} and~\eqref{eq:discreteLQOCP}. Most studies consider piecewise constant inputs and use numerical discretization methods to obtain the DT equivalent and derive the solution~\citep{ALT2016RegularizationandImplicitEuler,hager2000runge,cannon2000infinite,John2013EfficientImplementationofRiccati,sakizlis2005explicit}. 

On the other hand, LQ-OCPs have a strong connection with other advanced control algorithms, such as model predictive control (MPC). This relationship arises from using LQ-OCPs as the optimization problem within the linear MPC (LMPC) framework. In practice, nearly all LMPC applications employ a discrete-time LQ-OCP with diagonal weighting matrices, rather than the continuous-time formulation presented in~\eqref{eq:continuousLQOCP}. For the latter, the corresponding discrete weighting matrix $Q$ is full-element, as illustrated in~\eqref{eq:analyticQandq}. \cite{Karl1970IntroToStoContTheory,åström2011computer,franklin1990DigitalControl} explored the discrete equivalent of the LQ-OCP~\eqref{eq:continuousLQOCP} and described the analytic expressions of the desired equivalent, discrete weighting matrices $Q$ and $q_k$ described in~\eqref{eq:analyticQandq}. \cite{pannocchia2010computing,Pannocchia2015} proposed a novel computational procedure for the CT linear-quadratic regulator problem (CT-LQR). They pointed out that the discrete weighting matrix $Q$ can be computed with one exponential matrix instead of standard numerical discretization methods. The numerical experiments showed that the CT-LQR offers some advantages over the standard DT-LMPC. \cite{li2014event,wang2001continuous,morato2020model} also presented some other interesting results about CT-MPC. However, the existing literature has not investigated the case with time delays, which is common in real-world scenarios. Therefore, in this paper, we would like to investigate the DT implementation of CT-LMPC with time delays and compare it with the standard DT-LMPC. 

This paper is organized as follows. Section~\ref{sec:MPCDiscretization} introduces the control model and different objective functions of the CT-LMPC and present their discretization and implementation. Section~\ref{sec:NumericalExperiments} presents numerical experiments comparing and testing the proposed CT-LMPC and conventional DT-LMPC. The conclusions are given in Section~\ref{sec:NumericalExperiments}.

\section{Discretization of Model Predictive Control}
\label{sec:MPCDiscretization}
In this section, we introduce the design and discretization of CT-LMPC.
\subsection{Discretization of control model}
Transfer function models with time delays are common in industrial process, such as oil field, cement, etc. This paper considers the following MIMO transfer function model as the control model of CT-LMPC
\begin{subequations}
\label{eq:StochasticInputOutputModel}
\begin{align}
    \bs Z(s) &= G(s)U(s) + H(s) \bs W(s), \\
    \bs Y(s) &= \bs Z(s) + \bs V(s),
\end{align}    
\end{subequations}
with the deterministic and stochastic transfer functions $G(s)$ and $H(s)$
\begin{subequations}
\label{eq:SISOTransferFunctionModel}
\begin{align}
    & G(s) = \begin{bmatrix}
        g_{11}(s) & \cdots & g_{1 n_u}(s) \\ 
        \vdots & \ddots & \vdots \\ 
        g_{n_z 1}(s) & \cdots & g_{n_z n_u}(s)
    \end{bmatrix}, \, 
    g_{ij}(s) = \frac{b_{ij}(s)}{a_{ij}(s)}e^{-\tau_{ij}s}, 
    \\
    & H(s) = \begin{bmatrix}
        h_{11} (s)& \cdots & h_{1 n_w}(s) \\ 
        \vdots & \ddots & \vdots \\ 
        h_{n_z 1}(s) & \cdots & h_{n_z n_w}(s)
    \end{bmatrix}, \, 
    h_{ip}(s) = \frac{d_{ip}(s)}{c_{ip}(s)},
\end{align}    
\end{subequations}
where $U$ and $\bs W$ are the deterministic and stochastic inputs to the system. $\bs Z$ and $\bs Y$ are the output and measurement. $\bs V$ indicates the measurement noise. We assume that $g_{ij}(s)$ for $i=1,2,\ldots,n_z$, $j=1,2,\ldots,n_u$ with the input time delay $\set{\tau_{ij}} \in \mathbb{R}_0^+$ is proper, and $h_{ip}(s)$ for $p=1,2,\ldots, n_w$ is strictly proper and has no time delay.

\cite{HAGDRUP2016171} described a Noise-Separation (NS) state space representation. The input-output transfer function model~\eqref{eq:StochasticInputOutputModel} can be converted into a deterministic part $Z^d(s)$=$G(s)U(s)$ and a stochastic part $\bs Z^s(s)$=$H(s) \bs W(s)$
\begin{subequations}
\label{eq:ns-continuousStateSpace}
    \begin{align}
    & Z^d(s) \sim \left\{
    \begin{aligned}
     \dot x^d_{ij}(t) &= A_{c,ij}^d x_{ij}^d(t)+ B_{c,ij}^d u_j(t-\tau_{ij}), \\
        z^d_{ij}(t) &= C_{c,ij}^d x_{ij}^d(t) + D_{c,ij}^d u_j(t-\tau_{ij}),
     \end{aligned}
    \right.
    \label{eq:ns-DeterministicModel}
    \\    
    & \bs Z^s(s) \sim \left\{
    \begin{aligned}
     d \bs x^s(t) &= A_c^s \bs x^s(t) dt + B_c^s d \bs \omega(t), \\
        \bs z^s(t) &= C_c^s \bs x^s(t), \quad d\bs \omega(t) \sim N_{iid}(0, I dt),
     \end{aligned}
    \right.
    \label{eq:ns-StochasticModel}
    \end{align}
\end{subequations}
with the system variables 
\begin{subequations}
    \label{eq:ns-outputFunctions}
    \begin{align}
        & x^d = \begin{bmatrix}
        {x^d_{11}}' & {x^d_{21}}' & \ldots & {x^d_{n_z n_u}}'
        \end{bmatrix}', \quad  z^d_{i}(t) = \sum_{j=1}^{n_u} z^d_{ij}(t),
        \\
        & z^d = \begin{bmatrix}
            z^d_{1} & z^d_{2} & \ldots & z^d_{n_z}
        \end{bmatrix}',
         \qquad \bs z(t) = z^d(t) + \bs z^s(t), 
    \end{align}
\end{subequations}
where the system state $\bs x(t_0) = [x_0^d; \bs x_0^s] \sim N(\hat x_0, P_0)$ and $\bs x_0^s = 0$ and $P_0=$diag([$\bs 0$, $P^s_0$]). 

In~\eqref{eq:ns-continuousStateSpace}, we convert the deterministic part $Z^d(s)$ into $[i\times j]$ SISO time-delay state space models as they may have different time delays. Note that $\mathcal{G}_c =\set{A_{c,ij}^d, B_{c,ij}^d, C_{c,ij}^d,D_{c,ij}^d}$ and $\mathcal{H}_c=\set{A_{c}^s, B_{c}^s, C_{c}^s}$ are deterministic and stochastic system matrices.

Define the system sampling time $T_s$, the corresponding discrete-time systems are 
\begin{subequations}
\label{eq:ns-discreteStateSpace}
    \begin{alignat}{3}
    & Z^d(s) \sim \left\{
    \begin{aligned}
     & x^d_{k+1} = A^d x_k^d + B^du_k,\\
     & z^d_k = C^d x^d_k + D^d u_k, 
     \end{aligned}
    \right.
    \label{eq:ns-discreteStateSpace-determinsitic}
    \\
    & \bs Z^s(s) \sim \left\{
    \begin{aligned}
     & \bs x^s_{k+1} = A^s \bs x_k^s + \bs w_k,\\
     & \bs z^s_k = C^s \bs x^s_k,
     \end{aligned}
    \right.
    \label{eq:ns-discreteStateSpace-stochastic}
    \end{alignat}
\end{subequations}
and 
\begin{align}
    & \bs y_k = \bs z_k + \bs v_k, \quad 
    \begin{bmatrix}
        \bs w_k \\ \bs v_k
    \end{bmatrix} \sim N \left(\begin{bmatrix}
        0 \\ 0
    \end{bmatrix}, \begin{bmatrix}
        R_{ww} & 0 \\ 0 & R_{vv}
    \end{bmatrix} \right),
\label{eq:initialStateCovariance}
\end{align}    
where the deterministic system~\eqref{eq:ns-discreteStateSpace-determinsitic} are obtained by stacking the all SISO deterministic models~\eqref{eq:ns-DeterministicModel}. Note that $\mathcal{G}_d =\set{A^d, B^d, C^d,D^d}$ and $\mathcal{H}_d=\set{A^s, C^s}$ are corresponding discrete system matrices. The covariance $R_{ww} = \int_0^{T_s}e^{A_c^s t} B_c^s {B_c^s}'{e^{A_c^st}}'dt$ can be solved by the matrix exponential $\footnotesize{\begin{bmatrix}
    \Phi_{11} & \Phi_{12} \\ 0 & \Phi_{21}
\end{bmatrix}= exp(\begin{bmatrix}
    -A_c^s & B_c^s {B_c^s}' \\ 0 & {A_c^s}'
\end{bmatrix}})$ and $R_{ww}=\Phi_{22}'\Phi_{12}$ or other numerical methods~\citep{zhaz2023LQDiscretization,zhaz2023LQDiscretizationWithDelays}.

The main advantage of using the NS state space representation is that one may run the Kalman filter with \eqref{eq:ns-StochasticModel} and run the control algorithm with~\eqref{eq:ns-DeterministicModel}. 

Assume that the initial states of deterministic and stochastic models, i.e. $x^d$, $\hat x^s$ and $P^s$ are available. The filtering update of Kalman filter can be performed as 
\begin{subequations}
\label{eq:NS-KalmanFilter}
\begin{align}
    & \hat y^s_{k|k-1} = C^s \hat x^s_{k|k-1},
    && y_k^s = y_k - \hat z^d_k, \\
    &  e_k = y^s_{k} - \hat y^s_{k|k-1}, \quad && \hat x^s_{k|k} = \hat x^s_{k|k-1} + K_{fx} e_k,
\end{align}
with the measurement covariance and the Kalman gain 
\begin{equation}
       R_{e} = C^s P^s ({C^s})' + R_{ww}, \quad 
       K_{fx} = P^s ({C^s})' R_{e}^{-1},
\label{eq:KalmanFilterCoefficients}
\end{equation}
\end{subequations}
where $\hat z^d_k$ is estimated deterministic output obtained by~\eqref{eq:ns-discreteStateSpace-determinsitic}. The matrix $P^s$ indicates the stationary stochastic state error covariance obtained by the solution of discrete-time Algebraic Riccati Equation (DARE)
\begin{equation}
\label{eq:KalmanFileterDARE}
    P^s =  A^s P^s ({A^s})' - A^s P^s ({C^s})' R_e^{-1}C^s P^s ({A^s})' + R_{vv}.
\end{equation}     
The estimated output then can be computed as
\begin{subequations}
\begin{align}
    & \hat x_{k+j+1|k}^s = A^s \hat x_{k+j|k}^s, 
    &&  \hat{z}^s_{k+j|k} = C^s \hat{x}^s_{k+j|k},
    \\
    & \hat z_{k+j|k} = \hat{z}^d_{k+j|k} + \hat{z}^s_{k+j|k}, 
    && j = 1,\ldots, N.
\end{align}
\label{eq:ns-predictionProcess}
\end{subequations}
 
\subsection{Reference tracking and input regularization objectives}
Define the output and input reference tracking error $\tilde z(t)$ as 
\begin{equation}
\label{eq:trackingError}
    \tilde z(t) = \begin{bmatrix}
        z(t) - \bar z(t)\\ u(t) - \bar u(t)
    \end{bmatrix} = \begin{bmatrix}
        z(t) - \bar z_k\\ u(t) - \bar u_k
    \end{bmatrix}, \quad \text{for} \:\: t_k \leq t < t_{k+1},
\end{equation}
where we consider piecewise constant parameterization on references $\bar z(t) = \bar{z}_k$ and $\bar u(t) = \bar{u}_k$ for $t_k \leq t < t_{k+1}$. 

Note that the reference trajectory becomes $\bar z_{k+j} = \bar z(t) - \hat{z}^s_{k+j|k}$ for $t\in [t_{k+j}, t_{k+j+1})$ when using NS state space expressions.

We then can define the following CT LQ-OCP for minimizing the output and input reference tracking error as
\begin{subequations}
\label{eq:LQOCP-ReferenceTracking-stochastic}
\begin{alignat}{3}
    & \min_{x,u,z,\tilde z} \quad  && \phi_z + \phi_u=  E\set{\int_{t_0}^{t_0+T} l_{c\tilde z}(\tilde z(t)) dt} \\
    & s.t. && \bs x(t_0)\sim N(\hat x_0, P_0),~\eqref{eq:ZOHInput},~\eqref{eq:ns-continuousStateSpace}, ~\eqref{eq:ns-outputFunctions},~\eqref{eq:trackingError},
\end{alignat}
\end{subequations}
with the stage cost function $l_{\tilde z}(\tilde z(t))$
\begin{equation}
    l_{c\Tilde z}(\tilde z(t)) = \frac{1}{2} \norm{W_{\Tilde z} \Tilde z(t)}^2_2 = \frac{1}{2} \Tilde z (t)' Q_{c\tilde z} \Tilde z(t),
\end{equation}
where $Q_{c\tilde z} = \text{diag}\left( Q_{cz}, Q_{cu} \right) = W_{\tilde z}' W_{\tilde z}$. $Q_{cz} $ and $Q_{cu}$ are weighting matrices for the reference tracking and input regularization objectives.

The corresponding DT equivalent of~\eqref{eq:LQOCP-ReferenceTracking-stochastic} is
\begin{subequations}
\label{eq:DiscreteTimeLQOCP-ReferenceTracking} 
\begin{alignat}{5}
& \min_{x,u} \quad &&\phi_z + \phi_u = \sum_{k\in \mathcal{N}} l_{\Tilde z,k}(x_k,u_k)  \\
& s.t. && x_0 = \hat x_0,~\eqref{eq:ns-discreteStateSpace-determinsitic},~\eqref{eq:ns-predictionProcess}
\end{alignat}    
\end{subequations}
with the stage cost function $l_{\Tilde z, k} (x_k, u_k)$
\begin{equation}
    l_{\bar{z},k}(x_k,u_k) = \frac{1}{2} \begin{bmatrix} x_k \\ u_k \end{bmatrix}' Q \begin{bmatrix} x_k \\ u_k \end{bmatrix} + q_k' \begin{bmatrix} x_k \\ u_k \end{bmatrix} + \rho_k,
\label{eq:deterministic-stageCost}
\end{equation}
where  
\begin{equation}
    q_k =  M [\bar z_k; \bar u_k], \quad  \rho_k = l_{c\tilde z}( [\bar z_k; \bar u_k]) T_s.,  \quad k \in \mathcal{N}.
\label{eq:qkandrho_k}
\end{equation}

\cite{zhaz2023LQDiscretization,zhaz2023LQDiscretizationWithDelays} described the matrices of the DT system ($A$,$B$,$Q$,$M$) of~\eqref{eq:DiscreteTimeLQOCP-ReferenceTracking} as a differential equation system and introduced numerical methods to solve them. 



Note that $u=[u_0; u_1; \cdots; u_{N-1}]$ and $u_k = I_k u$ for $I_k = [0 \cdots  I  \cdots  0]$ and $k \in \mathcal{N}$. The system state $x_k$ can be expressed as 
\begin{equation}
     x_{k} = b_k + \Gamma_k u, \quad b_k = A^k x_0, \quad \Gamma_k = \sum_{i=0}^k A^{k-1-i} B I_i.
\label{eq:bkandGammak}
\end{equation}
We then rewrite the LQ-OCP~\eqref{eq:DiscreteTimeLQOCP-ReferenceTracking} in the form of QP
\label{eq:QP-ReferenceTracking} 
\begin{alignat}{5}
    \phi_z + \phi_u = \frac{1}{2} u'H_{\Tilde z} u + g_{\Tilde z}' u,
\end{alignat}
where 
\begin{equation}
    H_{\Tilde z} = \sum_{k=0}^{N-1} \begin{bmatrix}
        \Gamma_k \\ I_k
    \end{bmatrix}' Q \begin{bmatrix}
        \Gamma_k \\ I_k
    \end{bmatrix}, \: \:
    g_{\Tilde z} = \sum_{k=0}^{N-1} \begin{bmatrix}
        \Gamma_k \\ I_k
    \end{bmatrix}' \left( Q \begin{bmatrix}
        b_k \\ 0
    \end{bmatrix} + q_k \right).
\label{eq:Hzgz}
\end{equation} 

\subsection{Input ROM and economics objectives}
We then consider the LQ-OCP for penalizing the input rate-of-movement (ROM) and input cost (economics)
\begin{equation}
\label{eq:ContinuousTimeLQOCP-InputROM}
    \begin{split}
        \phi_{\Delta u} + \phi_{eco} &= \int_{t_0}^{t_0 + T} l_{c\Tilde u}(\dot u(t), u(t)) dt,
    \end{split}
\end{equation}
and the stage cost function is 
\begin{equation}
\label{eq:ContinuousTimeLQOCP-InputROM-stageCost}
    l_{c\Tilde u}(\dot u(t), u(t)) = \frac{1}{2} \norm{W_{c\Delta u} \dot u(t)}^2_2 + q_{ceco}' u(t),
\end{equation}
where $Q_{c\Delta u} = W_{c\Delta u}' W_{c\Delta u}$ and $q_{ceco}$ are the weight matrices of the input ROM and the economics objectives. 

The discretization of the input cost objective is intuitive if a zero-order hold (ZOH) parameterization is applied to the inputs. However, this is not the case for the input ROM objective function, as the input derivative is not even defined. \cite{Hagdrup2019MPCforSystems} described the discretization schemes on the input ROM penalty function using piecewise affine functions (FOH) and ZOH discretization. The author approved that the solution of the discrete problem is convergent toward the original problem's minimizer when taking the number of discretizations $N \rightarrow \infty$.

Lets consider the ZOH discretization on the input $u(t)=u_k$ for $t\in [t_k, t_{k+1})$ and the discrete approximation of the input ROM penalty is 
\begin{equation}
    \phi_{\Delta u} = \frac{1}{2}\int_{t_0}^{t_0 + T}\norm{\dot u(t)}^2_{Q_{c\Delta u}}dt = \frac{1}{2T_s} \sum_{k\in \mathcal N} \norm{u_k - u_{k-1}}^2_{Q_{c\Delta u}}.
\end{equation}
We then can perform the discretization of~\eqref{eq:ContinuousTimeLQOCP-InputROM} as
\begin{equation}
    \phi_{\Delta u} + \phi_{eco} = 
    \sum_{k \in \mathcal{N}} l_{\Tilde u, k} (u_{k}, u_{k-1}),
\end{equation}
with the stage cost function 
\begin{equation}
    \begin{split}
         l_{\Tilde u, k} (u_{k}, u_{k-1}) = \frac{1}{2} \begin{bmatrix}
            u_k \\ u_{k-1}
        \end{bmatrix}'
        \bar Q_{\Delta u}  \begin{bmatrix}
            u_k \\ u_{k-1}
        \end{bmatrix} + 
        q_{eco}' u_k,
    \end{split}
\end{equation}
where $u_{-1}$ indicates the input at $t_{0-1}$. The discrete-time weights $\bar Q_{\Delta u}$ and $q_{eco}$ are
\begin{equation}
    \bar Q_{\Delta u} = \begin{bmatrix}
        Q_{\Delta u} & - Q_{\Delta u} \\ 
        -Q_{\Delta u} & Q_{\Delta u}
    \end{bmatrix}, \: Q_{\Delta u} = \frac{Q_{c\Delta u}}{T_s}, \:
    q_{eco} = q_{ceco} T_s.
\end{equation}
Consequently, along with possible input box and input ROM constraints, we can express the input ROM and economics objective functions in the form of QP as 
\begin{subequations}
\label{eq:QP-InputROM}
 \begin{alignat}{3}
    & \min_{u} \quad && \phi_{\Delta u} + \phi_{u} = \frac{1}{2} u' H_{\Tilde u} u + g_{\Tilde u}' u
    \\
    & s.t. && u_{\min, k} \leq u_k \leq u_{\max, k}, \quad && k \in \mathcal{N}, \label{eq:inputConstraints}\\
& && \Delta u_{\min, k} \leq \Delta u_k \leq \Delta u_{\max, k}, \quad && k \in \mathcal{N}, \label{eq:inputROMConstraints}
\end{alignat}       
\end{subequations}
where the quadratic and linear term coefficients 
\begin{subequations}
\label{eq:Hdugdu}
\begin{align}
    & H_{\Tilde u} = \sum_{k=1}^{N-1} \begin{bmatrix}
        I_{k} \\ I_{k-1}
    \end{bmatrix}' \bar Q_{\Delta u} \begin{bmatrix}
        I_{k} \\ I_{k-1}
    \end{bmatrix} +  I_0' Q_{\Delta u}I_0,
    \label{eq:Hu}
    \\
    & g_{\Tilde u} =  \sum_{k=0}^{N-1} -I_0'  Q_{\Delta u} u_{-1} + I_k' q_{eco}. 
    \label{eq:gu}
\end{align}    
\end{subequations} 
\subsection{Soft output constraint penalty}
We then introduce the soft output constraints
\begin{subequations}
\label{eq:softConstraints}
\begin{align}
    & z_{k+j|k} \geq z_{\min, k+j|k}  - \xi_{k+j}, && k = 1,2,\ldots, N,  
    \\
    & z_{k+j|k} \leq z_{\max, k+j|k} + \eta_{k+j}, && k = 1,2,\ldots, N, 
    \\
    & \xi_{k+j} \geq 0, && k = 1,2,\ldots, N,  \label{eq:softConstraint01}
    \\
    & \eta_{k+j} \geq 0, && k = 1,2,\ldots, N, \label{eq:softConstraint02}
\end{align}    
where $\xi$ and $\eta$ are slackness variables. The output $z_k$ is subject to the deterministic system~\eqref{eq:ns-discreteStateSpace}. $z_{\min, k+j|k} = z_{\min, k+j} - \hat{z}^s_{k+j|k}$ and $z_{\max, k+j|k} = z_{\max, k+j} - \hat{z}^s_{k+j|k}$ are modified soft constraints. 
\end{subequations}

The corresponding penalty function
\begin{align}
    \phi_{\xi} + \phi_{\eta} = \int_{t_0}^{t_0+T} l_{c\xi}(\xi(t)) + l_{c\eta}(\eta(t)) dt,
\end{align}
with the stage cost functions 
\begin{align}
    & l_{c\xi}(\xi(t)) = \frac{1}{2} \norm{ W_{c\xi} \xi(t)}^2_2 + q_{c\xi}' \xi(t), 
    \\
    & l_{c\eta}(\eta(t)) = \frac{1}{2} \norm{ W_{c\eta} \eta(t)}^2_2 + q_{c\eta}' \eta(t), 
\end{align}
where $Q_{c\xi} = W_{c\xi}' W_{c\xi}$, $Q_{c\eta} = W_{c\eta}' W_{c \eta}$, $q_{c\xi}$ and $q_{c\eta}$ are weight matrices of the quadratic and linear penalty functions, respectively. We assume piece-wise constant $\xi(t)=\xi_k$ and $\eta(t)=\eta_k$ for $ t \in [t_k, \, t_k+T_s)$. The corresponding discrete equivalent is 
\begin{alignat}{3}
\label{eq:phist}
    \phi_{\xi} + \phi_{\eta}= \sum_{k=1}^{N} \frac{1}{2}(\norm{\xi_k}^2_{Q_{\xi}} + \norm{\eta_k}^2_{Q_{\eta}}) + q_{\xi}' \xi_k + q_{\eta}' \eta_k, 
\end{alignat}
where $Q_{\xi} = T_s Q_{c\xi}$, $q_{\xi} = T_s q_{c\xi}$, $Q_{\eta} = T_s Q_{c\eta}$ and $q_{\eta} = T_s q_{c\eta}$.

The DT penalty functions~\eqref{eq:phist} can be rewritten in the form of QP as
\begin{subequations}
\label{eq:QP-softConstraintPenalty}
    \begin{align}
        \phi_{\xi} + \phi_{\eta} = \frac{1}{2} \begin{bmatrix}
            \xi \\ \eta
        \end{bmatrix}' H_{\Tilde s}  \begin{bmatrix}
            \xi \\ \eta
        \end{bmatrix} + g_{\Tilde s}' \begin{bmatrix}
            \xi \\ \eta
        \end{bmatrix},
    \end{align}
\end{subequations}
where $\xi$=$[\xi_{1}; \xi_2; \cdots; \xi_{N}]$ and $\eta$=$[\eta_{1}; \eta_2; \cdots; \eta_{N}]$ and 
\begin{subequations}
 \begin{align}
    & H_{\xi} = \text{diag}({Q}_{\xi}; {Q}_{\xi}; \cdots; Q_{\xi}),
    \: g_{\xi} = \begin{bmatrix} q_{\xi} & q_{\xi} & \cdots & q_{\xi}
                 \end{bmatrix}',
    \\
    & H_{\eta} = \text{diag}({Q}_{\eta}; {Q}_{\eta}; \cdots; Q_{\eta}),      
    g_{\eta} = \begin{bmatrix}
                    q_{\eta} & q_{\eta} & \cdots & q_{\eta}
                \end{bmatrix}',
    \\
    & H_{\Tilde s} = \text{diag}({H}_{\xi}; {H}_{\eta}),
    \qquad \quad \:\:\: g_{\Tilde s} = \begin{bmatrix} g_{\xi} & g_{\eta} \end{bmatrix}'.
\end{align}   
\end{subequations}
\subsection{Design and implementation of CT-LMPC}
%
\begin{algorithm}[tb]
\caption{Design of CT-LMPC}
\label{algo:DesignOfLQMPC}
\begin{flushleft}
    \textbf{Require:} $(G(s),H(s),Q_{c\Tilde z},Q_{c\Delta u},Q_{c\xi}, Q_{c\eta}, q_{ceco}, q_{c\xi}, q_{c\eta},$ \\
    $R_{vv}, N, T_s)$ \\
\end{flushleft}
\begin{algorithmic}
\State Compute NS state space matrices $\mathcal{G}_c$, $\mathcal{H}_c$, $\mathcal{G}_d$ and $\mathcal{H}_d$ using~\eqref{eq:ns-continuousStateSpace} and~\eqref{eq:ns-discreteStateSpace}
\State Compute Kalman gain $K_{fx}$ using~\eqref{eq:KalmanFileterDARE}
\State Compute DT weight matrices $Q$ and $M$ using~\eqref{eq:deterministic-stageCost}
\State Compute $b=\set{b_k}$ and $\Gamma=\set{\Gamma_k}$ for $k\in \mathcal{N}$ using~\eqref{eq:bkandGammak}
\State Compute $H_{\Tilde z}$ of $\phi_{z}$ and $\phi_{u}$ using~\eqref{eq:Hzgz} 
\State Compute $H_{\Tilde u}$ of $\phi_{\Delta u}$ and $\phi_{eco}$ using~\eqref{eq:Hu}
\State Compute $H_{\Tilde s}$ and $g_{\Tilde s}$ of $\phi_{\xi}$ and $\phi_{\eta}$
using~\eqref{eq:QP-softConstraintPenalty}
\State Set $H=$ diag($[H_{\Tilde z}+H_{\Tilde u}; H_{\Tilde s}]$)
\end{algorithmic}
\begin{flushleft}
    \textbf{Return:} $(\mathcal{G}_d, \mathcal{H}_d, K_{fx}, Q, M, b, \Gamma, H, g_{\Tilde s})$ 
\end{flushleft}
\end{algorithm}
\begin{algorithm}[tb]
\caption{Implementation of CT-LMPC}
\label{algo:ImplementationOfLQMPC}
\begin{flushleft}
    \textbf{Require:} $(\bar{z},\bar{u}, u_{\min}, u_{\max}, \Delta u_{\min}, \Delta u_{max}, z_{\min}, z_{\max}$, \\ $u_{k-1}, y_{k}, \hat{x}^s_{k-1|k-1}, \mathcal{G}_d, \mathcal{H}_d, K_{fx}, Q, M, b, \Gamma, H, g_{\tilde s})$ 
\end{flushleft}
\begin{algorithmic}
\State Run Kalman filter for $\hat{x}^s_{k|k}$ and $\set{\hat z^s_{k+j|k}}$ using~\eqref{eq:ns-predictionProcess}
\For{\texttt{j=1,2,$\ldots$,N}}
    \State Update $\bar z_{k+j}$, $z_{\min,k+j}$ and $z_{\max,k+j}$ 
    with $\hat z^s_{k+j|k}$
    \State Set $q_{k+j} = M [\bar z_{k+j}; \bar u_{k+j}]$ 
\EndFor
\State Compute $g_{\Tilde z}$ using~\eqref{eq:Hzgz} 
\State Compute $g_{\Tilde u}$ using~\eqref{eq:gu}
\State Set $g = [g_{\tilde z}+g_{\tilde u}; g_{\tilde s}]$
\State Solve the QP~\eqref{eq:LQOCP-MPC} to get the optimal solution $u^*$
\State Let $u_{k} \leftarrow u^*(1:n_u)$
\end{algorithmic}
\begin{flushleft}
    \textbf{Return:} $(u_k, \hat{x}^s_{k|k})$ 
\end{flushleft}
\end{algorithm}
Combining the objective functions and constraints introduced previously, we have 
\begin{subequations}
\label{eq:LQOCP-MPC} 
\begin{alignat}{5}
& \min_{ \set{u,\xi,\eta}} \quad && \phi = \phi_{z} + \phi_u + \phi_{\Delta u} + \phi_{eco} + \phi_{\xi} + \phi_{\eta}
\\
& s.t. && x_0 = \hat x_0,~\eqref{eq:ns-discreteStateSpace-determinsitic},~\eqref{eq:ns-predictionProcess},~\eqref{eq:inputConstraints},~\eqref{eq:inputROMConstraints},~\eqref{eq:softConstraints},
\end{alignat}    
\end{subequations}
where the objective functions $\phi_{z}$, $\phi_u$, $\phi_{\Delta u}$, $\phi_{eco}$, $\phi_{\xi}$ and $\phi_{\eta}$ are the corresponding discrete equivalent of their original CT problems. 

The cost function $\phi$ can be expressed in the form of QP as 
\begin{equation}
    \phi = \frac{1}{2} \begin{bmatrix}
        u \\ \xi \\ \eta
    \end{bmatrix}'
     H  \begin{bmatrix}
        u \\ \xi \\ \eta
    \end{bmatrix} + g' \begin{bmatrix}
        u \\ \xi \\ \eta
    \end{bmatrix},
\label{eq:QP-LQMPC}
\end{equation}
where 
\begin{equation}
    H = \begin{bmatrix}
        H_{\Tilde z} + H_{\Tilde u} & 0 \\
        0 & H_{\Tilde s}
    \end{bmatrix}, \qquad 
    g = \begin{bmatrix}
        g_{\Tilde z} + g_{\Tilde u} 
        \\
        g_{\Tilde s}
    \end{bmatrix}.
\label{eq:QP-LQMPC-Coefficients}
\end{equation}
Consequently, we obtain the objective function~\eqref{eq:LQOCP-MPC} that is the discrete-time equivalent of the CT LQ-OCPs introduced in previous subsections. Algorithm~\ref{algo:DesignOfLQMPC} describes the design of CT-LMPC and the implementation of the CT-LMPC is illustrated by Algorithm~\ref{algo:ImplementationOfLQMPC}.

\section{Numerical Experiments}
\label{sec:NumericalExperiments}
\begin{figure*}[tb]
    \centering
    \begin{subfigure}[t]{0.33\textwidth} 
        \centering
        \includegraphics[width=1.\linewidth]{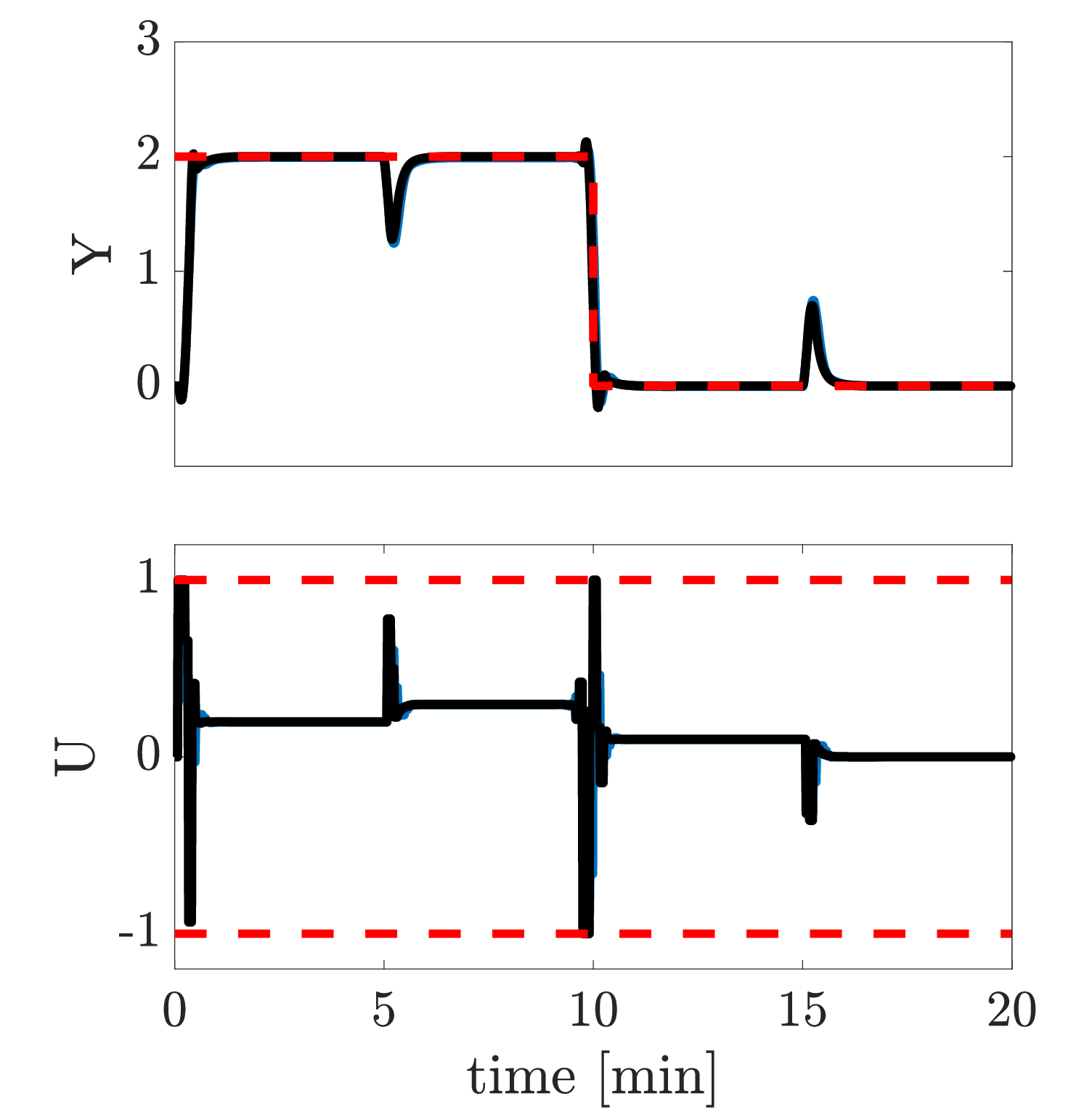}
        \caption{Deterministic sim. $T_s$ = 5 [s]}
        \label{fig:SISOExample-ClosedLoopSimulations-Det5Ts}
    \end{subfigure}%
    \hfill
    \begin{subfigure}[t]{0.33\textwidth}
        \centering
        \includegraphics[width=1.\linewidth]{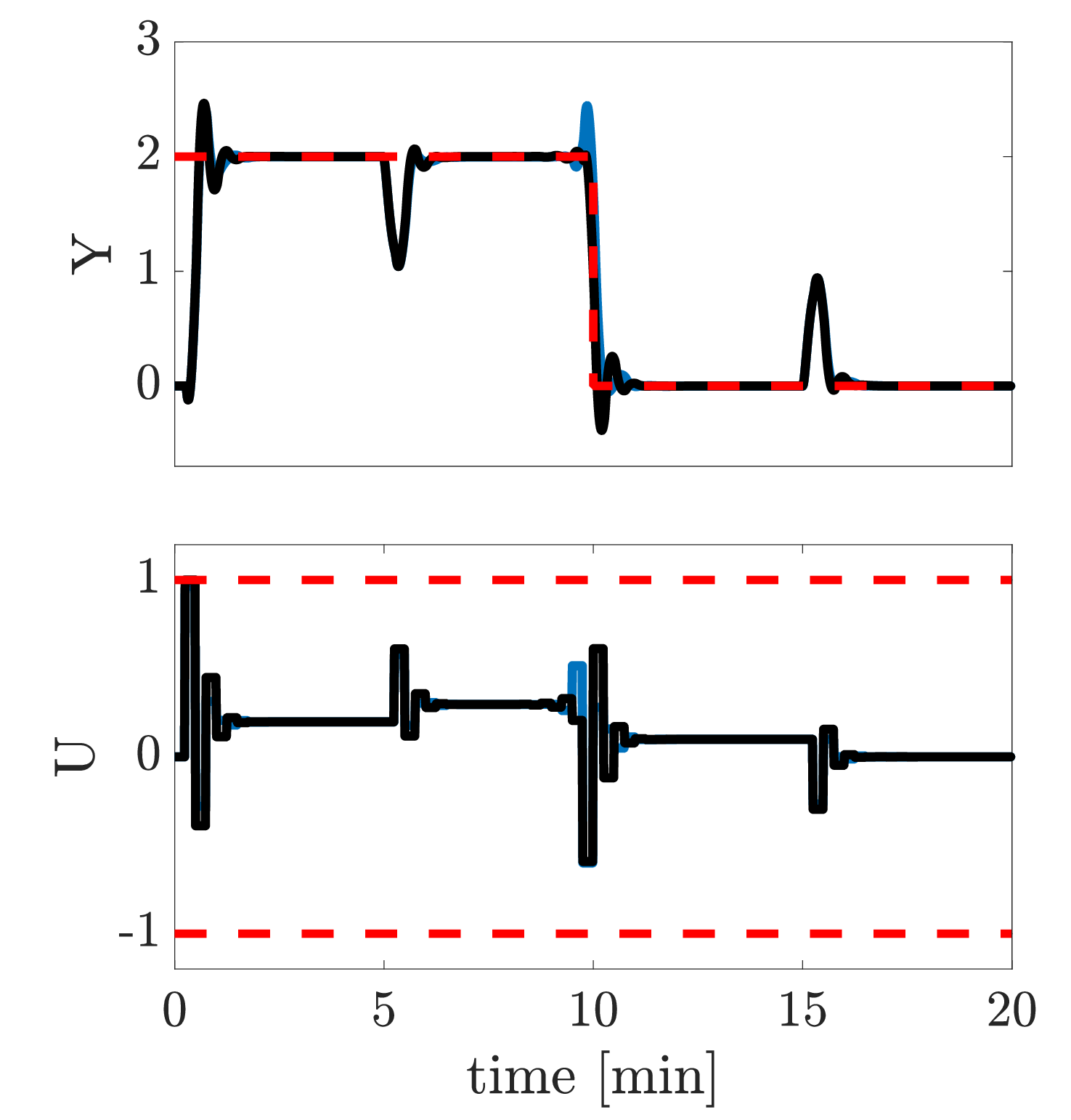}
         \caption{Deterministic sim. $T_s$ = 15 [s]}
         \label{fig:SISOExample-ClosedLoopSimulations-Det15Ts}
    \end{subfigure}%
    \hfill
    \begin{subfigure}[t]{0.33\textwidth} 
        \centering
        \includegraphics[width=1.\linewidth]{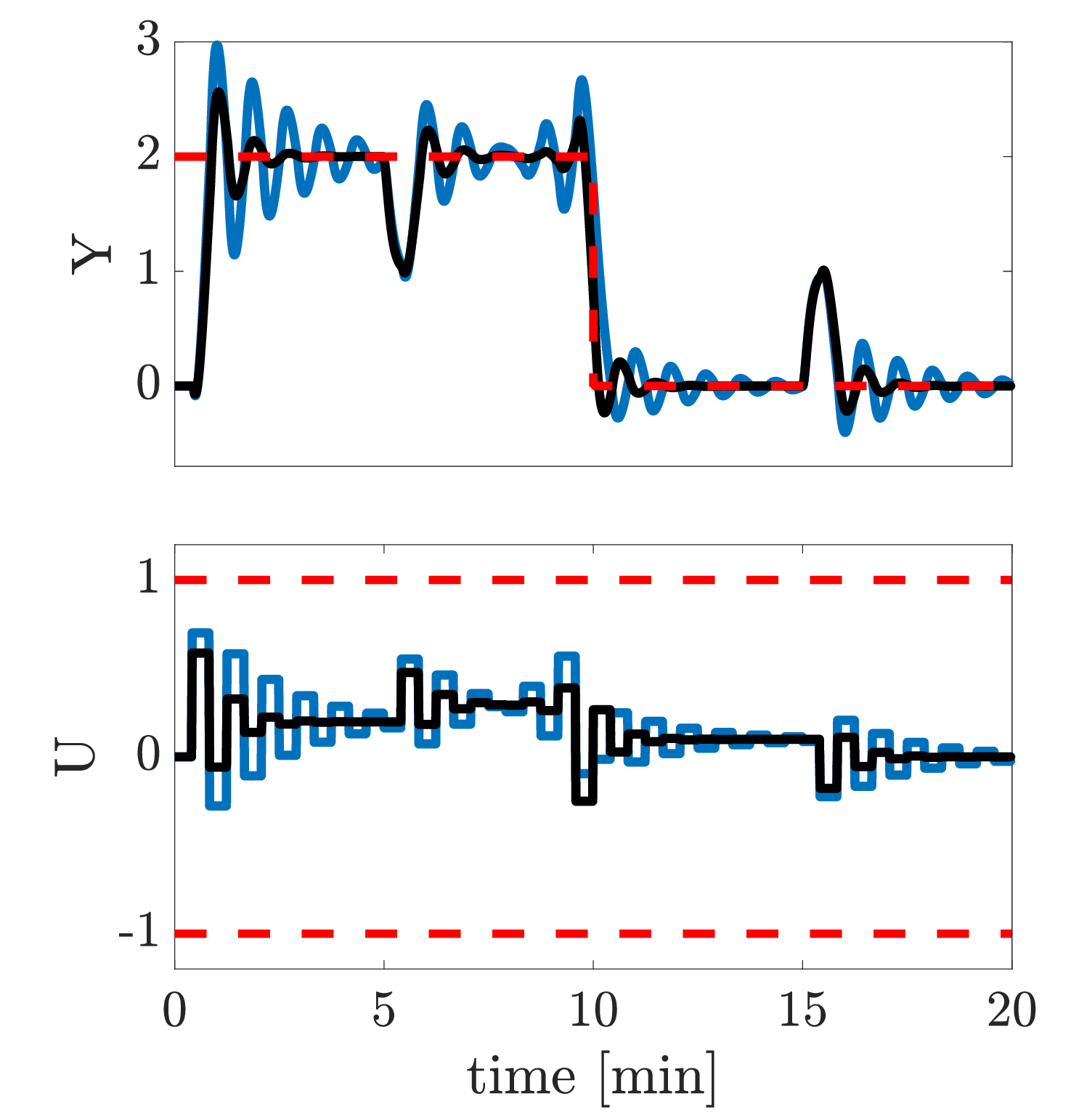}
        \caption{Deterministic sim. $T_s$ = 25 [s]}
        \label{fig:SISOExample-ClosedLoopSimulations-Det25Ts}
    \end{subfigure}%
    \hfill
    \begin{subfigure}[t]{0.33\textwidth} 
        \centering
        \includegraphics[width=1.\linewidth]{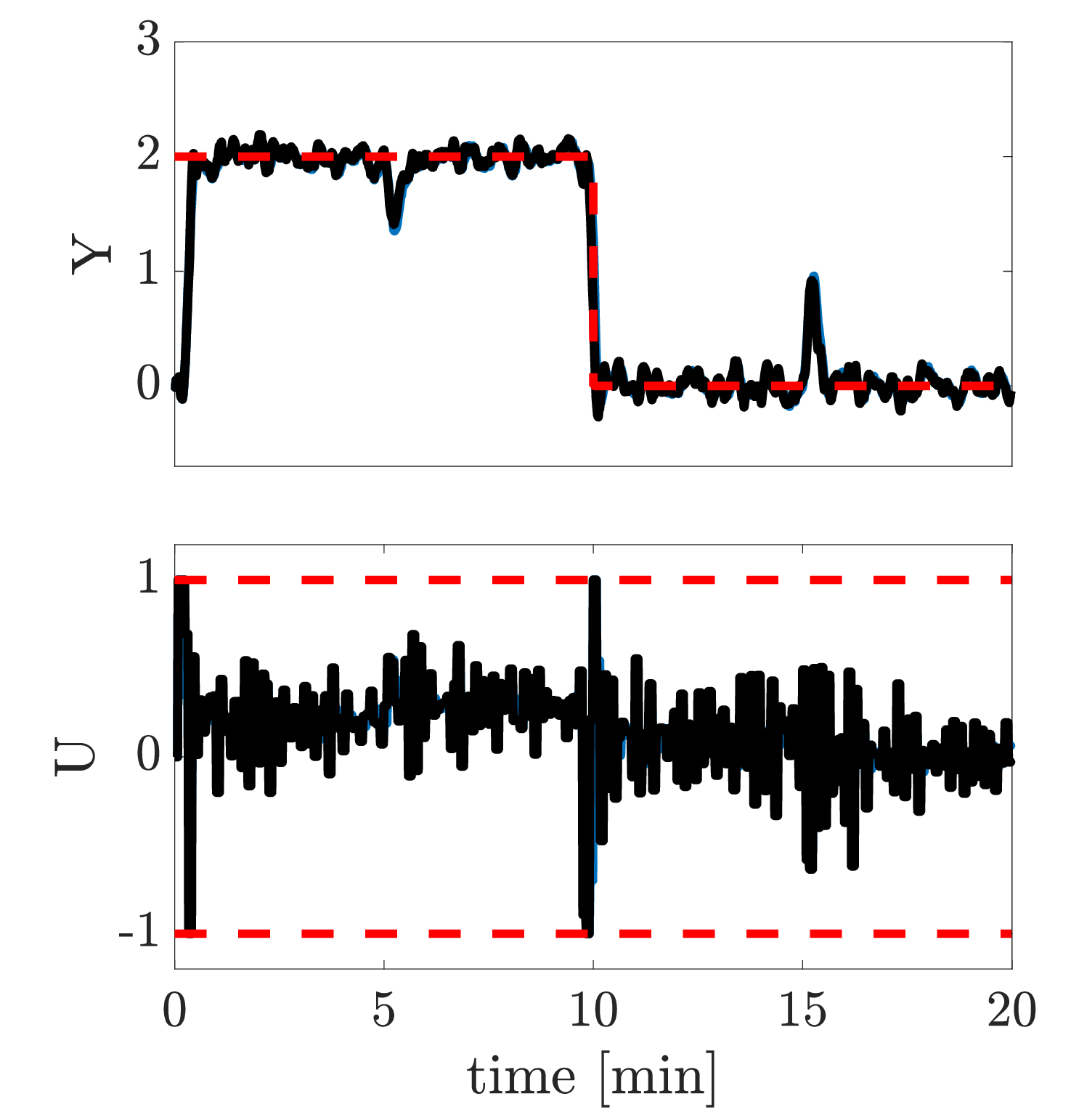}
        \caption{Stochastic sim. $T_s$ = 5 [s]}
        \label{fig:SISOExample-ClosedLoopSimulations-Stc5Ts}
    \end{subfigure}%
    \hfill
    \begin{subfigure}[t]{0.33\textwidth} 
        \centering
        \includegraphics[width=1.\linewidth]{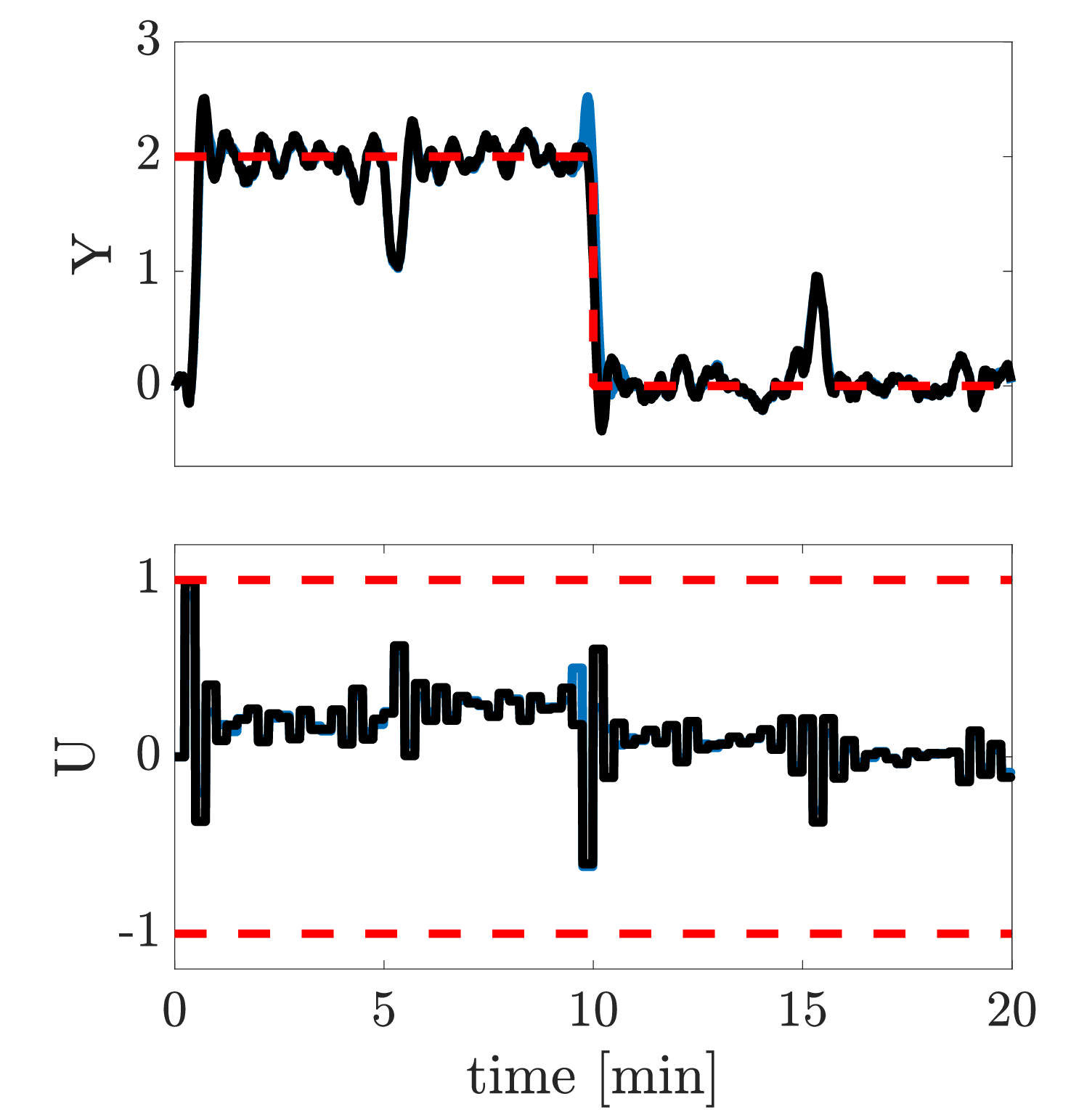}
        \caption{Stochastic sim. $T_s$ = 15 [s]}
        \label{fig:SISOExample-ClosedLoopSimulations-Stc15Ts}
    \end{subfigure}%
    \hfill
    \begin{subfigure}[t]{0.33\textwidth} 
        \centering
        \includegraphics[width=1.\linewidth]{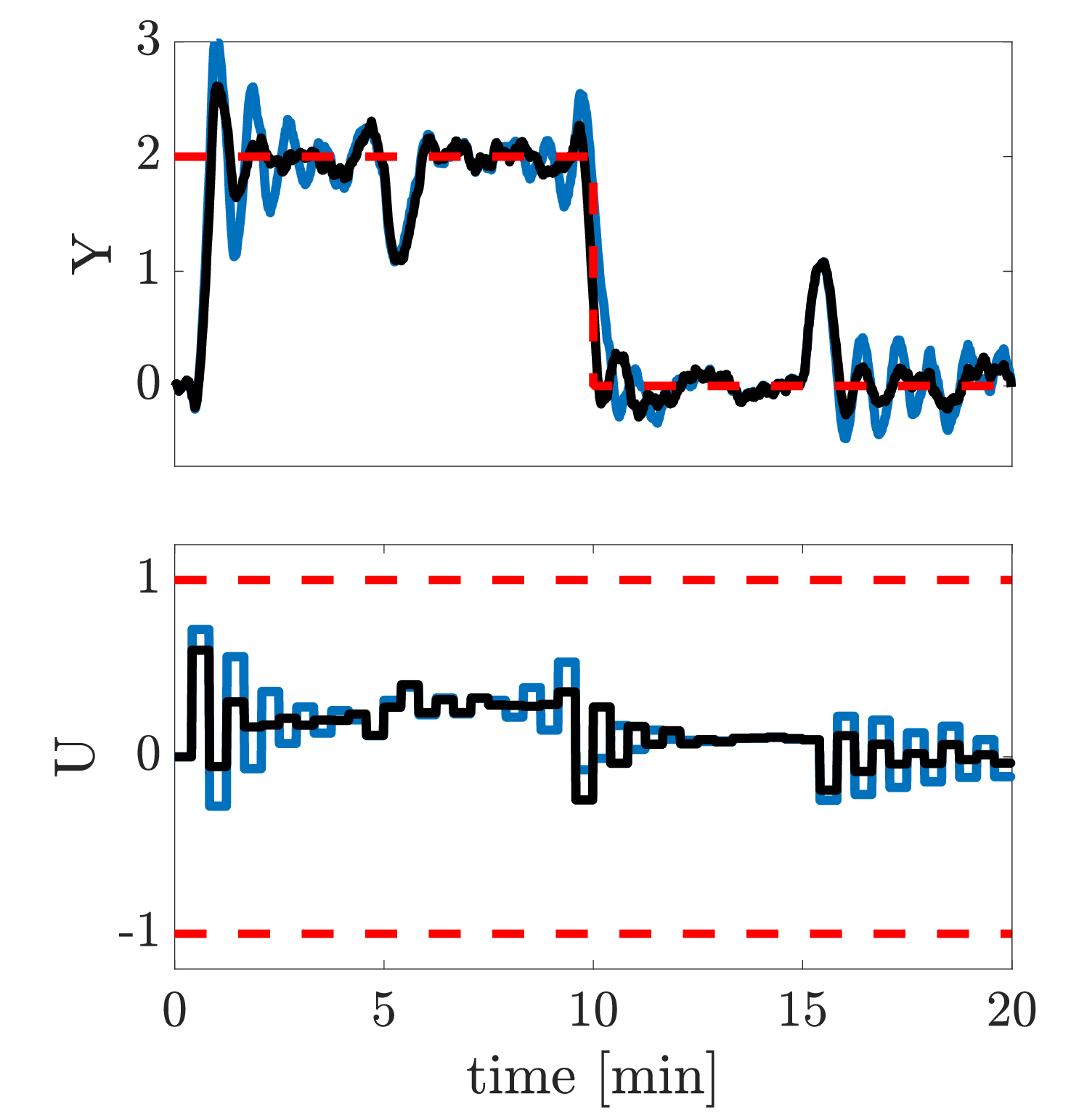}
        \caption{Stochastic sim. $T_s$ = 25 [s]}
        \label{fig:SISOExample-ClosedLoopSimulations-Stc25Ts}
    \end{subfigure}%
\caption{SISO example closed-loop simulations with different controller sampling times $T^c_{s} = 5, 15, 25$ [s]. The blue curves indicate the results obtained from DT-LMPC and the black curves are the results of CT-LMPC.}
\label{fig:SISOExample-ClosedLoopSimulations}
\end{figure*}
This section describes closed-loop simulations with the proposed CT-LMPC and compares it with the standard DT-LMPC. Consider the simulation model with the following expressions
\begin{subequations}
\label{eq:InputOutputModel}
\begin{align}
    & \bs Z(s) = G(s) U(s) + G_d(s) (D(s) + \bs W(s)), 
    \\
    & \bs Y(s) =  \bs Z(s) + \bs V(s),
\end{align}    
\end{subequations}
where $G_d(s)$ is the transfer function for the disturbance $D$ and uncertainty $\bs W$.

We discretize the above transfer function model with the sampling time $T_s$ as
\begin{subequations}
\label{eq:discreteTimeStateSpaceSimulator}
    \begin{align}
        \bs x_{k+1} & = A \bs x_k + B u_k + E d_k + G \bs w_k, \: \bs w_k \sim N_{iid}(0, R_{ww}), 
        \\
        \bs y_k & = C \bs x_k + D u_k + \bs v_k, \qquad \qquad \: \: \bs v_k \sim N_{iid}(0, R_{vv}),
    \end{align}
\end{subequations}
where the system state $x_k \sim N(\hat x_k, P)$ is available.

We develop the DT-LMPC based on the previous work by~\cite{HAGDRUP2016171} and the control model~\eqref{eq:StochasticInputOutputModel} for CT- and DT-LMPC is obtained by step-response experiments. It may be done using system parameters identification schemes introduced by~\cite{John2007PEM01,John2007PEM02,Daniel2013MPCTuning}. 
\subsection{Closed-loop simulation - a SISO example}
In the following, we perform a series of closed-loop simulations with CT- and DT-LMPC implementations for a SISO system. The SISO example considers the simulation model~\eqref{eq:InputOutputModel} with transfer functions
\begin{subequations}
\begin{align}
    & g(s) = \frac{10.12(-3.41s+1)e^{-2.5s}}{(15.9s+1)(24.2s+1)},
    \\
    & g_d(s) = \frac{-0.5}{(5.8s+1)(4.7s+1)}.
\end{align}
\end{subequations}
The transfer function model is converted into the state space~\eqref{eq:discreteTimeStateSpaceSimulator} with $T_s=1$ [s]. The unknown disturbance $d_k$ = 2.0 for $5 \leq t \leq 15$ [min]. The covariances are $R_{ww}=R_{vv}=0.02^2$.

The estimated control model~\eqref{eq:StochasticInputOutputModel} for the two MPCs is 
\begin{equation}
    \hat{g}(s) = \frac{10.12(-3.58s+1)e^{-2.5s}}{(18.9s+1)(22.2s+1)}, \: \:
    \hat h(s) = \frac{1}{s}\frac{0.6}{(s+1)},
\end{equation}
where we select $\hat h(s)$ as an integrator to ensure the controller can handle step-type error~\citep{MUSKE2002617,Pannocchia2003Offsetfree}. 
\begin{figure*}[bt]
    \centering
    \begin{subfigure}[t]{0.50\textwidth} 
        \centering
        \includegraphics[width=1.\linewidth]{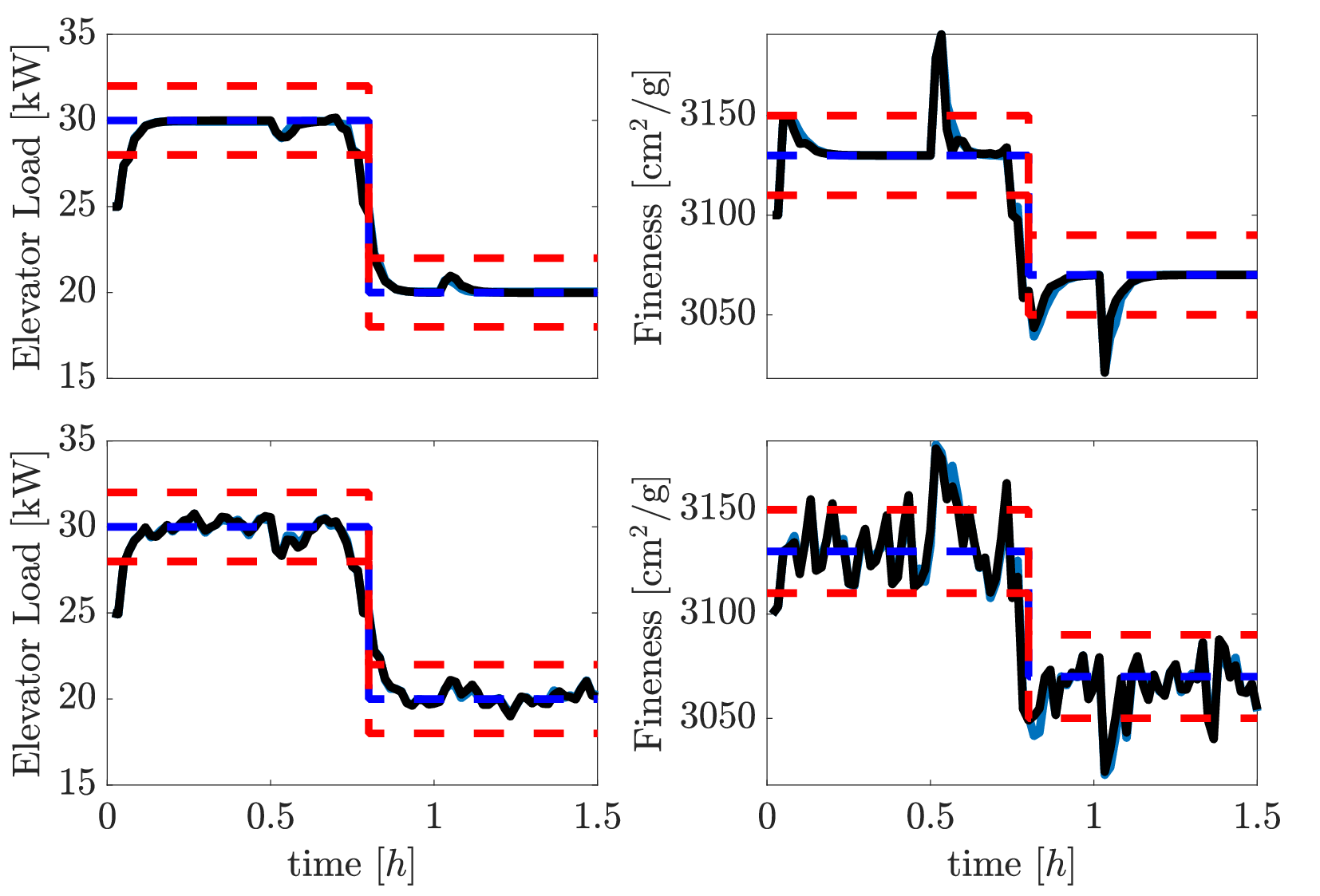}
        \caption{The outputs of deterministic and stochastic sim.}
        \label{fig:MIMOExample-ClosedLoopSimulations-Outputs}
    \end{subfigure}%
    \hfill
    \begin{subfigure}[t]{0.50\textwidth} 
        \centering
        \includegraphics[width=1.\linewidth]{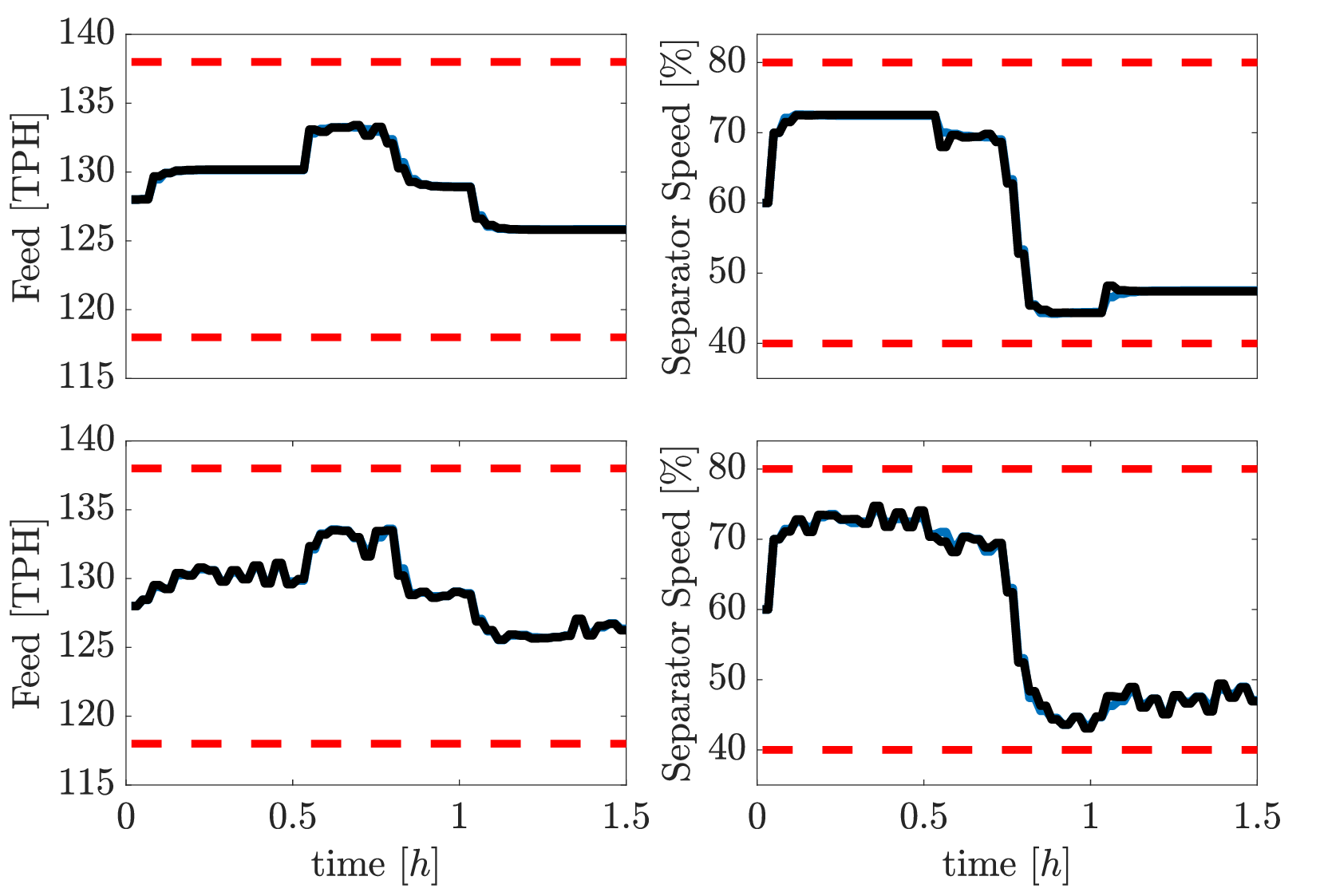}
        \caption{The inputs of deterministic and stochastic sim.}
        \label{fig:MIMOExample-ClosedLoopSimulations-Inputs}
    \end{subfigure}%
\caption{The closed-loop simulations of a simulated cement mill system with DT-LMPC (blue curves) and CT-LMPC (black curves) implementations. The simulator and controller sampling times are $T_s$ = 1 [min] and $T_s^c$ = 2 [min]}
\label{fig:MIMOExample-ClosedLoopSimulations}
\end{figure*}
In the SISO example, we consider the reference tracking objective $\phi_z$ and input ROM objective $\phi_{\Delta u}$ for two predictive controllers. 
The cost functions $\phi_z$ and $\phi_{\Delta u}$ for CT-LMPC are obtained based on Algorithm~\ref{algo:DesignOfLQMPC}. For DT-LMPC, $\phi_z$ and $\phi_{\Delta u}$ are designed as
\begin{equation}
\label{eq:dtMPC-ReferenceTrackingInputROM}
    \phi = \sum_{k=0}^{N-1} \norm{z_{k+1} - \bar{z}_{k+1}}^2_{Q_{cz}} + \norm{u_k - u_{k-1}}^2_{Q_{c\Delta u}},
\end{equation}
The prediction and control horizon is $N=20$. The output target is $\bar{z} = 2$ for $t \leq 7.5$ [min] and $\bar z = -2$ for the remaining time. The weight matrices are $Q_{cz} = 20$ and $Q_{c\Delta u} = 1$, and the input constraints are $u_{\min}$ = -1 and $u_{\max}$ = 1. 

To demonstrate the contrast between CT- and DT-LMPC, we discretize both using different controller sampling times $T^c_{s}$ = 5, 15, 25 [s]. Fig.~\ref{fig:SISOExample-ClosedLoopSimulations} shows closed-loop simulations of the SISO system with CT- and DT-LMPC implementations. In Fig.~\ref{fig:SISOExample-ClosedLoopSimulations-Det5Ts} and~\ref{fig:SISOExample-ClosedLoopSimulations-Stc5Ts}, both CT-LMPC (black curves) and DT-LMPC (blue curves) effectively control the SISO system, allowing the outputs to follow the given target (red dashed line). Two overshoots occur at $t=5$ and $15$ [min] due to an unknown disturbance, and the two MPCs can handle it after a few iterations. In this case, the simulation results show little difference between the two MPCs. Upon increasing $T^c_s$ to 15 [s], additional overshoots emerge around 1 and 10 [min] in Fig.~\ref{fig:SISOExample-ClosedLoopSimulations-Det15Ts} and~\ref{fig:SISOExample-ClosedLoopSimulations-Stc15Ts}. Notably, the second overshoot at $t=10$ [min] for DT-LMPC is larger than that for CT-LMPC. Moreover, the overshoots resulting from the unknown disturbance are more significant compared to the previous scenario, yet both MPCs can manage it. Fig.~\ref{fig:SISOExample-ClosedLoopSimulations-Det25Ts} and~\ref{fig:SISOExample-ClosedLoopSimulations-Stc25Ts} show the closed-loop simulations when $T_c^s$ = 25 [s]. Here, we can see a significant decrease in the closed-loop performance of both MPCs compared to the previous two cases. Both CT- and DT-LMPC exhibit oscillations and struggle to follow the target. However, the amplitude of the output from CT-LMPC is considerably smaller than that from DT-LMPC. CT-LMPC maintains control of the output to the desired target after 1-2 iterations, while DT-LMPC fails to do so.

The SISO example simulation results show that a large sampling time can degrade the closed-loop performance of the predictive controller and even lead to system instability. For the same system parameters, the proposed CT-LMPC has a better closed-loop performance than the traditional DT-LMPC, and this gap will increase with the increase of the sampling time.

\subsection{Closed loop simulation - a MIMO example}
The MIMO example concerns the cement mill system introduced by~\cite{PRASATH2010ApplicationofSoftConstrainedMPC,Daniel2013MPCTuning}, where the cement mill system is described by transfer functions 
\begin{subequations}
 \begin{align}
    & G(s) = \begin{bmatrix}
        \frac{0.62e^{-5s}}{(45s+1)(8s+1)} & \frac{0.29(8s+1)e^{-1.5s}}{(2s+1)(38s+1)} \\ 
        \frac{-15e^{-5s} }{60s+1}& \frac{5e^{-0.1s}}{(14s+1)(s+1)}
    \end{bmatrix},
    \\
    & G_d(s) = \begin{bmatrix}
        \frac{-1.0e^{-3s}}{(32s+1)(21s+1)} \\ \frac{60}{(30s+1)(20s+1)}
    \end{bmatrix}.
\label{eq:CementMillModel}
\end{align}   
\end{subequations}
The system inputs $u$ = [feed flow rate (TPH); separator speed (\%)] and the system outputs $z$ = [elevator load (kW); fineness (cm$^2$/g)]. The system disturbance $D$ is the clinker hardness (HGI). The transfer functions of the estimated control model~\eqref{eq:StochasticInputOutputModel} for the two MPCs are
\begin{subequations}
    \begin{align}
        & \hat G(s) = \begin{bmatrix}
        \frac{0.8e^{-5s}}{(30s+1)(15s+1)} & \frac{0.45e^{-2s}}{30.0s+1} \\ 
        \frac{-17.7e^{-5s} }{(65s+1)(15s+1)}& \frac{9.4 e^{-0.3s} }{15s+1}
    \end{bmatrix},
    \\
    & \hat H(s) = \text{diag} \left( \left[ \frac{1}{s} \frac{0.5}{s+1}; \frac{1}{s} \frac{1}{s+1} \right] \right),
    \end{align}
\end{subequations}
where we manually introduce the plant-model mismatch on the deterministic model $\hat G(s)$.

The cement mill system is converted into a DT state space model~\eqref{eq:discreteTimeStateSpaceSimulator} with the sampling time $T_s$ = 1 [min]. The simulation time $T_{sim}=1.5$ [h]. The related covariances are selected as $R_{ww}$ = 1.0 and $R_{vv}$ = diag([0.1; 50]). The unknown disturbance $d_k$ = 8 for $t \in [0.5, 1]$ [h] and $d_k$ = 0 for the rest of the time. 

For the CT- and DT-LMPC, the controller sampling time $T_s^c=2$ [min], the prediction and control horizon is $N=60$. The weight matrices for the reference tracking, input ROM, economics and soft output constraints penalty are $Q_{cz}=$ diag($200, 10$), $Q_{c\Delta u}=$ diag($20, 10$), $q_{ceco}=[2;1]$, $Q_{c\xi}=Q_{c\eta}=$ diag($2000, 100$) and $q_{c\xi}=q_{c\eta}=[20; 1]$, respectively. The input hard constraints are $u_{\min}=[-10; -20]$, $u_{\max}=[10, 20]$, $\Delta u_{\min}=[-5, -10]$ and $\Delta u_{\max}=[5, 10]$. The soft output constraints $z_{\min}=[-2, -20]$ and $z_{\max}=[2, 20]$ and are deviation variables. 

Fig.~\ref{fig:MIMOExample-ClosedLoopSimulations} illustrates the closed-loop simulations of a simulated cement mill system with CT- and DT-LMPC implementations. In Fig.~\ref{fig:MIMOExample-ClosedLoopSimulations-Outputs}, both CT-LMPC (black curves) and DT-LMPC (blue curves) can regulate the outputs towards the desired targets (blue dashed lines) and within the given bounds for most of the time. An overshoot, attributed
to the plant-model mismatch in $y_2$ (Fineness), is observed. The two predictive controllers can handle the problem and bring the output back to the target after a few iterations. 
At $t=0.5$ and $1.5$ hours, two overshoots on the output occur due to the unknown disturbance. Both CT- and DT-LMPC can reject the unknown disturbance. Additionally, an overshoot arises at around t = 0.75 [h] on the Fineness due to a step change in the output targets. DT-LMPC exhibits a slightly larger overshoot than CT-LMPC.

The MIMO example closed-loop simulations show that the proposed CT-LMPC and DT-LMPC exhibit highly similar closed-loop performance. They even converge to identical optimal solutions when using appropriate sampling time and tuning parameters. This outcome aligns with the earlier SISO example presented in Fig.~\ref{fig:SISOExample-ClosedLoopSimulations-Det5Ts} and~\ref{fig:SISOExample-ClosedLoopSimulations-Stc5Ts}, reinforcing the consistency of our findings.

\section{Conclusions}
\label{sec:Conclusions}
This paper presents the design, discretization and implementation of CT-LMPC. We introduce different objective functions of CT-LMPC and show their discretization and implementation. From the numerical experiments, we notice that the proposed CT-LMPC has the following features: 
\begin{itemize}
    \item [1.] When using an appropriate sampling time and well-tuned weight matrices, the proposed CT-LMPC has a highly similar control performance as the traditional DT-LMPC, and they will converge to an identical optimal solution
    \item [2.] The proposed CT-LMPC has a better closed-loop performance than the traditional DT-LMPC with the same system parameters when the controller sampling time is large. This gap will increase with the increase in controller sampling time
\end{itemize}
Besides, compared with conventional DT-LMPC, it is more reasonable to design a continuous-time LQ-OCP for MPC and then discretize it, as this is how we do it in NMPC.

\bibliography{ref/reference}                         

\end{document}